\documentclass{article}
\usepackage{amsfonts}

\newtheorem{theorem}{Theorem}

\newtheorem{corollary}[theorem]{Corollary}

\newtheorem{remark}[theorem]{Remark}

\input{tcilatex}
\begin{document}

\begin{center}
{\Large Identities on the }${\large k}${\Large -ary Lyndon words related to
a family of zeta functions}

Irem Kucukoglu$^{1,a}$ and Yilmaz Simsek$^{1,b}$

$^{a}$ikucukoglu@akdeniz.edu.tr

$^{b}$ysimsek@akdeniz.edu.tr

$^{1}$Department of Mathematics, Faculty of Science University of Akdeniz
TR-07058 Antalya, Turkey

\bigskip%

\textbf{Abstract}
\end{center}

The main aim of this paper is to investigate and introduce relations between
the numbers of $k$-ary Lyndon words and unified zeta-type functions which
was defined by Ozden et al.\cite[p. 2785, Definition 3]%
{OzdenSimsekSrivastava2010}. Finally, we give some identities on\ generating
functions for the numbers of $k$-ary Lyndon words and some special numbers
and polynomials such as the Apostol-Bernoulli numbers and polynomials,
Frobenius-Euler numbers, Euler numbers and Bernoulli numbers.



\section{Introduction}

Throughout this paper, we consider the number of $k$-ary Lyndon words of
length $n$, $L_{k}\left( n\right) $ as follows \cite{ThomasBOOK}:

\begin{equation}
L_{k}\left( n\right) =\frac{1}{n}\sum_{d|n}\mu \left( \frac{n}{d}\right)
k^{d},  \label{Witt}
\end{equation}%
where the arithmetic function $\mu $ is the M\"{o}bius function defined as
follows \cite{APOSTOLbook}:

\[
\mu \left( 1\right) =1; 
\]%
Let $p_{1}$, $p_{1}$,$\ldots $, $p_{k}$ be $k$ distinct primes. If $n>1$,
write $n=p_{1}^{a_{1}}p_{2}^{a_{2}}\ldots p_{k}^{a_{k}}$. Then%
\begin{eqnarray*}
\mu \left( n\right) &=&\left( -1\right) ^{k}\text{ if }a_{1}=a_{2}=\ldots
=a_{k}=1\text{, } \\
&&\text{that is, if }n\text{ is the product of }k\text{ distinct primes,} \\
\mu \left( n\right) &=&0\text{ if }n\text{ is the product of non distinct
primes.}
\end{eqnarray*}

In \cite{Lyndon}, Lyndon words are studied as standard lexicographic
sequences. According to \cite[p. 36]{ThomasBOOK}, \textit{a }$k$\textit{-ary
necklace is an equivalence class of }$k$\textit{-ary strings under rotation.
As a representative of such an equivalence class which is taken the smallest
in the lexicographical order. A period }$n$\textit{\ necklace representative
with }$n$\textit{\ digits is called a Lyndon word}.

In addition to counting $k$-ary Lyndon words of length $n$, Equation (\ref%
{Witt}) is well known as Witt's formula which is used to count the number of
monic irreducible polynomials of degree $n$ over Galois field (\textit{cf. }%
\cite{Buchanan}) and As we mentioned previously in \cite{Kucukoglu1}, it is
know that there are really interesting connections between this formula and
dimension formula for the homogeneous subspaces of the free Lie algebra (%
\textit{cf}. \cite{Buchanan}, \cite{Kang}, \cite{Petrogradsky}) and the rank
of the free abelian quotient (\textit{cf.} \cite{Witt1937}, \cite{Lyndon}).
Furthermore, it is also called necklace polynomial (\textit{cf.} \cite%
{Metropolis1}). For further information about $L_{k}\left( n\right) $ and
table including numerical values of $L_{k}\left( n\right) $, the reader can
consult \cite{Kucukoglu1}, its references and also the references cited in
each of these earlier works.

In \cite{Kucukoglu1}, the authors gave the following explicit formula for
the numbers of $L_{k}\left( p\right) $ for $p$ is a prime number, $m\in 
\mathbb{N}
$:%
\begin{equation}
L_{k}\left( p^{m}\right) =\frac{k^{p^{m-1}}\left( k^{p^{m-1}\left( {p-1}%
\right) }-1\right) }{p^{m}}.  \label{Formula}
\end{equation}

In \cite{Kucukoglu1}, the authors defined ordinary generating functions for
the numbers of $k$-ary Lyndon words of prime length $p$, $L_{k}\left(
p\right) $ for prime $p$, as follows:

\textit{Let }$p$\textit{\ is a prime number }and $m=1$ in the special case
of Equation (\ref{Formula})%
\begin{equation}
\mathfrak{f}_{L}(t,p)=\sum_{k=2}^{\infty }L_{k}\left( p\right) t^{k}=\frac{1%
}{p}\sum_{k=2}^{\infty }\left( k^{p}-k\right) t^{k},  \label{GF1}
\end{equation}%
(\textit{cf.} \cite{Kucukoglu1}). However, we modify (\ref{GF1}) as follows:

\[
f_{L_{y}}(t,p)=\sum_{k=0}^{\infty }L_{k}\left( p\right) t^{k}.
\]

\section{Relation between $f_{L_{y}}(t,p)$ and a family of zeta functions}

In this section, our aim is to give some identities on the generating
functions for the numbers of $k$-ary Lyndon words of length prime related a
family of zeta-type function, the Apostol-Bernoulli numbers and polynomials,
Frobenius-Euler numbers and Euler numbers.

The Apostol-Bernoulli numbers, $\mathcal{B}_{k}\left( z\right) $ are defined
by means of the following generating functions:%
\[
\frac{t}{ze^{t}-1}=\sum\limits_{k=0}^{\infty }\frac{\mathcal{B}_{k}\left(
z\right) }{k!}t^{k},
\]%
(\textit{cf. \cite{ApostolLERCHzeta}}). The Apostol Benoulli polynomials, $%
\mathcal{B}_{k}\left( x,z\right) $ are also defined by means of the
following generating functions:%
\[
\frac{te^{tx}}{ze^{t}-1}=\sum\limits_{k=0}^{\infty }\frac{\mathcal{B}%
_{k}\left( x,z\right) }{k!}t^{k},
\]%
where\ $\left\vert t\right\vert <2\pi $ when $z=1$ and $\left\vert
t\right\vert <\left\vert \log z\right\vert $ when $z\neq 1$ and $z\in 
\mathbb{C}$. From this generating function, one can easily see that%
\[
\mathcal{B}_{m}\left( x,z\right) =\sum\limits_{j=0}^{m}\left( 
\begin{array}{c}
m \\ 
j%
\end{array}%
\right) x_{j}^{m-j}\mathcal{B}_{j}\left( z\right) ,
\]%
and%
\[
\mathcal{B}_{m}\left( 0,z\right) =\mathcal{B}_{m}\left( z\right) .
\]

By using the above generating function, several of Apostol Bernoulli numbers
and polynomials are given as follows, respectively (\textit{cf. \cite%
{ApostolLERCHzeta}}):%
\begin{eqnarray*}
\mathcal{B}_{0}\left( z\right)  &=&0, \\
\mathcal{B}_{1}\left( z\right)  &=&\frac{1}{z-1}, \\
\mathcal{B}_{2}\left( z\right)  &=&\frac{-2z}{\left( z-1\right) ^{2}}, \\
\mathcal{B}_{3}\left( z\right)  &=&\frac{3z\left( z+1\right) }{\left(
z-1\right) ^{3}}, \\
\mathcal{B}_{4}\left( z\right)  &=&\frac{-4z\left( z^{2}+4z+1\right) }{%
\left( z-1\right) ^{4}}, \\
\mathcal{B}_{5}\left( z\right)  &=&\frac{5z\left( z^{3}+11z^{2}+11z+1\right) 
}{\left( z-1\right) ^{5}}, \\
\mathcal{B}_{6}\left( z\right)  &=&\frac{-6z\left(
z^{4}+26z^{3}+66z^{2}+26z+1\right) }{\left( z-1\right) ^{6}}
\end{eqnarray*}%
and%
\begin{eqnarray*}
\mathcal{B}_{0}\left( z,x\right)  &=&0, \\
\mathcal{B}_{1}\left( z,x\right)  &=&\frac{1}{z-1}, \\
\mathcal{B}_{2}\left( z,x\right)  &=&\frac{1}{z-1}x-\frac{2z}{\left(
z-1\right) ^{2}}, \\
\mathcal{B}_{3}\left( z,x\right)  &=&\frac{9z\left( z-1\right) }{\left(
z-1\right) ^{3}}x^{2}-\frac{6z}{\left( z-1\right) ^{2}}x+\frac{3z\left(
z-1\right) }{\left( z-1\right) ^{3}},
\end{eqnarray*}%
(\textit{cf. \cite{ApostolLERCHzeta}, \cite{HuMin}, \cite%
{OzdenSimsekSrivastava2010}, \cite{Ozden2014}, \cite{SimsekAMC2010}, \cite%
{SrivastavaChoi1}, \cite{SrivastavaOzden2012}, \cite{SrivastavaChoi2}}).

In the above numerical computation of the Apostol-Bernoulli numbers, we
observe that all of this numbers are rational functions of parameter $z$. $%
z=1$ is a pole of these functions.

The following generating function of the unification of the Bernoulli, Euler
and Genocchi polynomials, $\mathcal{Y}_{n,\beta }\left( x;k,a,b\right) $,
which was recently defined by Ozden \cite{Ozden2010} for $k\in \mathbb{N}%
_{0}:=\mathbb{N}\cup \{0\}(\mathbb{N}:=\{1,2,3,\ldots \})$; $a,b\in \mathbb{R%
}^{+}$; $\beta \in \mathbb{C}$:%
\begin{equation}
\frac{2^{1-k}t^{k}e^{tx}}{\beta ^{b}e^{t}-a^{b}}=\sum\limits_{n=0}^{\infty }%
\mathcal{Y}_{n,\beta }\left( x;k,a,b\right) \frac{t^{n}}{n!}
\end{equation}%
where $\left\vert t+b\log \left( \frac{\beta }{a}\right) \right\vert <2\pi $%
; $x\in \mathbb{R}$ and note that%
\[
\mathcal{Y}_{n,\beta }\left( k,a,b\right) =\mathcal{Y}_{n,\beta }\left(
0;k,a,b\right) =\mathcal{Y}_{n,\beta }\left( 1;k,a,b\right) .
\]

The following equation of the unified zeta-type functions $\zeta _{\beta
}\left( s,x;k,a,b\right) $, which was recently defined by Ozden et al.\cite[%
p. 2785, Definition 3]{OzdenSimsekSrivastava2010}:%
\begin{equation}
\zeta _{\beta }\left( s,x;k,a,b\right) =\left( \frac{-1}{2}\right)
^{k-1}\sum\limits_{n=0}^{\infty }\frac{\beta ^{bn}}{a^{b\left( n+1\right)
}\left( n+x\right) ^{s}}  \label{zeta1}
\end{equation}%
where $\beta ,s\in \mathbb{C}$ with $Re\left( s\right) <1$ and $|\beta |<1$
and observe that if\ $x=1$, then%
\[
\zeta _{\beta }(s:k,a,b)=\zeta _{\beta }\left( s,1;k,a,b\right) =\left( 
\frac{-1}{2}\right) ^{k-1}\sum\limits_{n=1}^{\infty }\frac{\beta ^{b\left(
n-1\right) }}{a^{b\left( n+1\right) }n^{s}}.
\]

Ozden et al.\cite[p. 2789, Theorem 7]{OzdenSimsekSrivastava2010} also proved
the following relation for $n\in \mathbb{N}$ and $k\in \mathbb{N}_{0}$:%
\begin{equation}
\zeta _{\beta }(1-n:x;a,b)=\left( -1\right) ^{k}\frac{\left( n-1\right) !}{%
\left( n+k-1\right) !}\mathcal{Y}_{n+k-1,\beta }\left( x;k,a,b\right) .
\label{zetaYR}
\end{equation}

\begin{remark}
Setting $s\rightarrow -m$, $\beta \rightarrow t$, $x=0$, $k=a=b=1$ in (\ref%
{zeta1}), we have%
\begin{equation}
\zeta _{\beta }(-m,0;1,1,1)=\sum\limits_{n=0}^{\infty }t^{n}n^{s}
\label{ZetaSumm}
\end{equation}%
and also setting $1-n\rightarrow -m$, $\beta \rightarrow t,$ $x=0$, $k=a=b=1$
in (\ref{zetaYR}), we get%
\begin{equation}
\zeta _{t}(-m:0;1,1)=-\frac{\mathcal{Y}_{1+m,t}\left( 0;1,1,1\right) }{1+m}
\label{ZetaYSp}
\end{equation}
\end{remark}

It is well-know that $\mathcal{Y}_{n,\beta }\left( x;1,1,1\right) $ reduce
to the Apostol-Bernoulli polynomials, $\mathcal{B}_{n}\left( x,\beta \right) 
$ and Apostol-Bernoulli numbers $\mathcal{B}_{n}\left( \beta \right) $,
respectively. Thus 
\begin{equation}
\mathcal{Y}_{n,\beta }\left( 0;1,1,1\right) =\mathcal{B}_{n}\left( 0,\beta
\right) =\mathcal{B}_{n}\left( \beta \right) .  \label{YApostol}
\end{equation}

Hurwitz-Lerch zeta function $\Phi \left( z,s,a\right) $ is defined by (%
\textit{cf}. \cite[p. 121 et seq.]{SrivastavaChoi2}, \cite[p. 194 et seq.]%
{SrivastavaChoi2}):%
\[
\Phi \left( z,s,a\right) =\sum\limits_{n=0}^{\infty }\frac{z^{n}}{\left(
n+a\right) ^{s}},
\]%
where $a\in \mathbb{C}\backslash \mathbb{Z}_{0}^{-}$; s $\in \mathbb{C}$
when $\left\vert z\right\vert <1$; $Re\left( s\right) >1$ when $\left\vert
z\right\vert =1$. One can also easily see that a relation between $\zeta
_{\beta }\left( s,x;k,a,b\right) $ and $\Phi \left( z,s,a\right) $ is given
by%
\[
\zeta _{\beta }\left( s,x;k,a,b\right) =\frac{\left( \frac{-1}{2}\right)
^{k-1}}{a^{b}}\Phi \left( \frac{\beta ^{b}}{a^{b}},s,x\right) 
\]%
(\textit{cf. }\cite[Eq-(4.14)]{OzdenSimsekSrivastava2010}). Hurwitz--Lerch
zeta function is related to not only Riemann zeta function and the Hurwitz
zeta function (see, for details, \cite[Chapter 2]{SrivastavaChoi2}, see also 
\cite{SrivastavaOzden2012}, \cite{SimsekAMC2010}):%
\[
\zeta \left( s\right) =\sum\limits_{n=1}^{\infty }\frac{1}{n^{s}}=\Phi
\left( 1,s,1\right) ,
\]%
and 
\[
\zeta \left( s,a\right) =\sum\limits_{n=0}^{\infty }\frac{1}{\left(
n+a\right) ^{s}}=\Phi \left( 1,s,a\right) .
\]

\begin{remark}
Let $n\geq 1$. Then in \cite{KimJKMS2008}, \cite{SimsekFPT2013} and \cite%
{SimsekASCM2016}, one can see that%
\begin{equation}
\mathcal{B}_{n}\left( z\right) =\frac{n}{z^{-1}}\mathcal{H}_{n-1}\left( 
\frac{1}{z}\right) ,  \label{FrobeniusApostol}
\end{equation}%
where $\mathcal{H}_{n}\left( z\right) $ denotes the Frobenius-Euler numbers
which are defined by means of the following generating function:%
\[
\frac{1-z}{e^{t}-z}=\sum\limits_{n=0}^{\infty }\mathcal{H}_{n}\left(
z\right) \frac{t^{n}}{n!},
\]%
for $z=-1$, we have $\mathcal{H}_{n}\left( -1\right) =E_{n}$ which is
defined by%
\[
\frac{2}{e^{t}+1}=\sum\limits_{n=0}^{\infty }E_{n}\left( z\right) \frac{t^{n}%
}{n!},
\]%
(cf. \cite{JangPak2002}-\cite{SrivastavaChoi2}; and the references cited
therein).
\end{remark}

Now, by combining Equation (\ref{ZetaSumm}) with Equation (\ref{ZetaYSp})
and Equation (\ref{YApostol}), we obtain the following explicit formula of
generating functions for the numbers of $k$-ary Lyndon words in terms of the
Apostol-Bernoulli numbers:

\begin{theorem}
\label{fLYApostol}\textit{Let }$p$\textit{\ is a prime number. Then}%
\[
f_{L_{y}}(t,p)=\frac{\mathcal{B}_{2}\left( t\right) }{2p}-\frac{\mathcal{B}%
_{p+1}\left( t\right) }{p\left( p+1\right) }
\]%
where $\mathcal{B}_{p+1}\left( z\right) $ denotes the Apostol-Bernoulli
numbers.
\end{theorem}

\begin{remark}
If we substitute $p=2$ into Theorem \ref{fLYApostol}, we arrive at 
\[
f_{L_{y}}(t,3)=\frac{t^{2}}{\left( 1-t\right) ^{3}}
\]%
which was given in \cite[p. 3]{Kucukoglu1}.
\end{remark}

\begin{remark}
Observe that degree of polynomial in the numerator of $f_{L_{y}}(t,2)$ is
lower than its the denominator. Also, when $p=3$ into Theorem \ref%
{fLYApostol}, we also arrive at 
\[
f_{L_{y}}(t,3)=\frac{2t^{2}}{\left( t-1\right) ^{4}}.
\]
\end{remark}

By combining Theorem \ref{fLYApostol} with Equation (\ref{FrobeniusApostol}%
), we give the following Remark:

\begin{theorem}
\textit{Let }$p$\textit{\ is a prime number. Then}%
\begin{equation}
f_{L_{y}}(t,p)=\frac{t\left( \mathcal{H}_{1}\left( \frac{1}{t}\right) -%
\mathcal{H}_{p}\left( \frac{1}{t}\right) \right) }{p}  \label{FrobeLy}
\end{equation}%
where $\mathcal{H}_{p}\left( \frac{1}{t}\right) $ denotes the
Frobenius-Euler number.
\end{theorem}

Since $\mathcal{H}_{n}\left( -1\right) =E_{n}$, then Equeation (\ref{FrobeLy}%
) is reduced to the following Corollary:

\begin{corollary}
\[
f_{L_{y}}(-1,p)=\frac{E_{p}-E_{1}}{p}
\]%
where $E_{n}$ denotes Euler numbers.
\end{corollary}

\section{Further identities related to Bernoulli numbers}

In this section, we give some identities on\ the numbers of $k$-ary Lyndon
words related to the Apostol-Bernoulli numbers. Now, we recall definition
the Bernoulli polynomials $B_{n}\left( x\right) $ which are defined by means
of the following generating function:%
\[
\frac{te^{tx}}{e^{t}-1}=\sum\limits_{n=0}^{\infty }B_{n}\left( x\right) 
\frac{t^{n}}{n!},
\]%
where $\left\vert t\right\vert <2\pi $ and also%
\[
B_{n}=B_{n}\left( 0\right) 
\]%
which denotes the Bernoulli numbers (\textit{cf.} \cite{JangPak2002}-\cite%
{SrivastavaChoi2}; and the references cited therein).

The sum of the powers of integers is related to the Bernoulli numbers and
polynomials:%
\begin{equation}
\sum_{k=0}^{m}k^{n}=\frac{B_{n+1}\left( m+1\right) -B_{n+1}}{n+1}
\label{BernoulliR}
\end{equation}%
(\textit{cf.} \cite{DjordjevicMilovanovic2014}, \cite{SrivastavaChoi2}, \cite%
{SrivastavaManochaTreatise1984}; see also the references cited in each of
these earlier works).

After applying mobius inversion formula to Equation (\ref{GF1}), we have%
\begin{equation}
k^{n}=\sum_{d|n}dL_{k}\left( d\right) .  \label{AfterMobius}
\end{equation}

Hence, summing Equation (\ref{AfterMobius}) over all $0\leq k\leq m$ and
combining with Equation (\ref{BernoulliR}), we obtain the following Theorem:

\begin{theorem}
Let $n\geq 1$. Then%
\[
\sum_{d|n}d\sum_{k=0}^{m}L_{k}\left( d\right) =\frac{B_{n+1}\left(
m+1\right) -B_{n+1}}{n+1},
\]%
where $B_{n+1}\left( m\right) $ denotes Bernoulli polynomials.
\end{theorem}

\textbf{Acknowledgments. }The present investigation was supported by\textit{%
\ Scientific Research Project Administration of Akdeniz University}.


\begin{thebibliography}{99}
\bibitem{ApostolLERCHzeta} T. M. Apostol, On the Lerch zeta function,\ 
\textit{Pacific Journal of Mathematics}\textbf{\ 1 (}1951), pp. 161-167.

\bibitem{APOSTOLbook} T. M. Apostol, \textit{Introduction to Analytic Number
Theory, }Narosa Publishing, Springer Verlag, New Delhi, Chennai, Mumbai,
1998.

\bibitem{Buchanan} H.L. Buchanan, A. Knopfmacher, M.E. Mays, On the
cyclotomic identity and related product expansions, \textit{Australas. J.
Combin.} \textbf{8} (1993), pp. 233-245.

\bibitem{ThomasBOOK} T. W. Cusick and P. Stanica, \textit{Cryptographic
Boolean Functions and Applications}, London: Academic Press, Elsevier, 2009.

\bibitem{DjordjevicMilovanovic2014} G. B. Djordjevic and G. V. Milovanovic,
Special classes of polynomials, University of Nis, Faculty of Technology
Leskovac, 2014.

\bibitem{HuMin} S. Hu, M.-S. Kim, Two closed forms for the Apostol-Bernoulli
polynomials, arXiv:1509.04190.

\bibitem{JangPak2002} \textbf{L. C. Jang and H. K. Pak, Non-archimedean
integration associated with }$q$\textbf{-Bernoulli numbers, Proc. Jangjeon
Math. Soc. 5(2) (2002), pp. 125-129.}

\bibitem{Kang} S.-J. Kang, M.-H. Kim. 1996. Free Lie algebras, generalized
Witt formula, and the denominator identity. \textit{J. Algebra} \textbf{183}
(1996), no. 2, pp. 560-594.

\bibitem{KimJKMS2008} \textbf{T. Kim, S.-H. Rim, Y. Simsek, and D Kim, On
the analogs of Bernoulli and Euler numbers, related identities and zeta and }%
$l$\textbf{-functions, J. Korean Math. Soc. 45(2) (2008), pp. 435-453.}

\bibitem{Kucukoglu1} I. Kucukoglu, Y. Simsek, On k-ary Lyndon Words And
Their Generating Functions, to appear in \textit{AIP Conf. Proc. of ICNAAM}%
(2016).

\bibitem{LuoSrivastava2005} \textbf{Q.-M. Luo, H. M. Srivastava, Some
generalizations of the Apostol-Bernoulli and Apostol-Euler polynomials,\ J.
Math. Anal. Appl. 308 (2005), pp. 290-302.}

\bibitem{Lyndon} R. Lyndon, On Burnside problem I, Trans.\textit{\ American
Math. Soc.} \textbf{77} (1954), pp. 202-215.

\bibitem{Metropolis1} N. Metropolis, G.-C. Rota, Witt Vectors and the
Algebra of Necklaces, \textit{Adv. in Math.} \textbf{50} (1983), pp. 95-125.

\bibitem{Ozden2010} H. Ozden, Unification of generating function of the
Bernoulli, Euler and Genocchi numbers and polynomials, \textit{Amer. Inst.
Phys. Conf. Proc.} \textbf{1281} (2010), 1125-1128.

\bibitem{OzdenSimsekSrivastava2010} H. Ozden, Y. Simsek, H.M. Srivastava, A
unified presentation of the generating functions of the generalized
Bernoulli, Euler and Genocchi polynomials, \textit{Comput. Math. Appl.} 
\textbf{60} (2010), pp. 2779-2787.

\bibitem{Ozden2014} H. Ozden, Y. Simsek, Unified presentation of p-adic
L-functions associated with unification of the special numbers, \textit{Acta
Math. Hungar.} \textbf{144 }(\textbf{2}) (2014), pp. 515-529.

\bibitem{OzdenSimsek2014Apostol} \textbf{H. Ozden, Y. Simsek, Modification
And Unification Of The Apostol-Type Numbers And Polynomials And Their
Applications, \textit{Appl. Math. Comput.} 235 (2014), pp. 338-351.}

\bibitem{Petrogradsky} V.M. Petrogradsky, Witt's formula for restricted Lie
algebras, \textit{Adv. Appl. Math.} \textbf{30} (2003), pp. 219-227.

\bibitem{SimsekKimJP2004} \textbf{Y. Simsek, T. Kim, D.W. Park, Y.S. Ro,
L.C. Jang, S. Rim, An Explicit Formula For The Multiple Frobenius-Euler
Numbers And Polynomials, \textit{JP J. Algebra Number Theory Appl.} (2004),
no.3, pp. 519-529.}

\bibitem{SimsekqAnalog2005} \textbf{Y. Simsek, q-Analogue of the twisted
l-Series and q-Twisted Euler Numbers, \textit{J. Number Theory} 100(2)
(2005), pp. 267-278.}

\bibitem{SimsekASCM2007} \textbf{Y. Simsek, O. Yurekli, V. Kurt, On
interpolation functions of the twisted generalized Frobenius--Euler numbers, 
\textit{Adv. Stud. Contemp. Math.} 15 (2007), no. 2, pp. 187-194.}

\bibitem{SimsekAMC2010} Y. Simsek and H. M. Srivastava, A family of p-adic
twisted interpolation functions associated with the modified Bernoulli
numbers, \textit{Appl. Math. Comput.} \textbf{216} (2010), pp. 2976-2987.

\bibitem{SimsekFPT2013} \textbf{Y. Simsek, Generating functions for
generalized Stirling type numbers, Array type polynomials, Eulerian type
polynomials and their applications, \textit{Fixed Point Theory A.} 2013,
2013:87}.

\bibitem{SimsekASCM2016} \textbf{Y. Simsek, Apostol type Dahee numbers and
Polynomials,\ to appear in Adv. Stud. Contemp. Math. 26(3) (2016).}

\bibitem{SrivastavaManochaTreatise1984} H. M. Srivastava and H. L. Manocha, 
\textit{A Treatise on Generating Functions}, Ellis Horwood Limited
Publisher, Chichester, 1984.

\bibitem{SrivastavaChoi1} H.M. Srivastava, J. Choi, \textit{Series
Associated with the Zeta and Related Functions}, Kluwer Academic Publishers,
Dordrecht, Boston and London, 2001.

\bibitem{SrivastavaOzden2012} H. M. Srivastava, H. Ozden, I. N. Cangul and
Y. Simsek, A unified presentation of certain meromorphic functions related
to the families of the partial zeta type functions and the L-functions, 
\textit{Appl. Math. Comput.} \textbf{219 }(2012), pp. 3903-3913.

\bibitem{SrivastavaChoi2} H.M. Srivastava, J. Choi,\textit{\ Zeta and q-Zeta
Functions and Associated Series and Integrals}, Elsevier Science Publishers,
Amsterdam, London and New York, 2012.

\bibitem{Witt1937} E. Witt, Treue Darstellung Liescher Ringe, \textit{J.
Reine Angew. Math.} \textbf{177 }(1937), pp.152-160.
\end{thebibliography}
\end{document}